\documentclass[a4paper,11pt,reqno]{amsart}
\usepackage{color}
\usepackage[hyphens]{url}
\usepackage[round,comma,sort&compress]{natbib}
\usepackage[cmfonts]{iba}
\usepackage{enumerate}
\usepackage{colortbl}
\usepackage{booktabs}
\usepackage{ifpdf}
\usepackage{graphics}
\newcommand\vealpha{{\boldsymbol{\alpha}}}
\newcommand\velambda{{\boldsymbol{\lambda}}}
\newcommand\vehatlambda{{\boldsymbol{\hat\lambda}}}
\newcommand\vetildelambda{{\boldsymbol{\tilde\lambda}}}

\DeclareMathOperator{\ind}{ind}
\DeclareMathOperator{\td}{td}
\DeclareMathOperator{\cone}{cone}
\newcommand\C{\mathbf C}

\usepackage[breaklinks=true,colorlinks,citecolor=blue]{hyperref}

\title[A primal Barvinok algorithm]{A primal Barvinok algorithm\\
  based on irrational decompositions}
\author{Matthias~K\"oppe}
\address{Otto-von-Guericke-Universit\"at Magdeburg, Department of
  Mathematics, Institute for Mathematical Optimization (IMO),
  Universit\"atsplatz 2, 
  39106 Magdeburg, Germany} 
\email{mkoeppe@imo.math.uni-magdeburg.de}
\thanks{This work was supported by a Feodor Lynen Research Fellowship from the
  Alexander von Humboldt Foundation.  The author also acknowledges the
  hospitality of Jes\'us De Loera and the Department of Mathematics of the
  University of California, Davis.
  The author is grateful to David Einstein for pointing out a mistake in
  an earlier version of this paper.}
\date{$\relax$Revision: 1.41 $ - \ $Date: 2006/09/25 11:54:39 $ $}
\begin{document}

\begin{abstract}
  We introduce variants of Barvinok's algorithm for counting lattice points in
  polyhedra.  The new algorithms are based on irrational signed decomposition
  in the primal space and the construction of rational generating functions
  for cones with low index.  We give computational results that show that the
  new algorithms are faster than the existing algorithms by a large factor. 
\end{abstract}

\maketitle

\section{Introduction}

\textsc{Twelve years} have passed since Alexander Barvinok's 
amazing algorithm for counting lattice points in polyhedra was published
\citep{bar}.  In 
the mean time, efficient implementations \citep{latte1,
  verdoolaege-et-al:counting-parametric} were designed, which helped to make
Barvinok's algorithm a practical tool in many applications in discrete
mathematics.  The implications of Barvinok's technique, of course, reach far
beyond the domain of combinatorial counting problems: For example,
\cite{deloera-etal-2005:five-ip-algorithms} pointed out applications in Integer
Linear Programming, and
\cite{deloera-hemmecke-koeppe-weismantel:mixedintpoly-fixeddim,deloera-hemmecke-koeppe-weismantel:intpoly-fixeddim}
obtained a fully polynomial-time approximation scheme ({\small FPTAS}) for
optimizing arbitrary polynomial functions over 
the mixed-integer points in polytopes of fixed dimension.

Barvinok's algorithm first triangulates the supporting cones of all vertices
of a polytope, to obtain simplicial cones.  Then, the simplicial cones 
are recursively decomposed into unimodular cones.  It is essential that one
uses \emph{signed decompositions} here; triangulating these cones is not good
enough to give a polynomiality result.  The rational generating functions of
the resulting unimodular cones can then be written down easily.  Adding and
subtracting them according to the inclusion-exclusion principle and the theorem
of \citet{Brion88} gives the rational generating function of the polytope.
The number of lattice points in the polytope can finally be obtained by
applying residue techniques on the rational generating function.

The algorithm in the original paper \citep{bar} worked explicitly with 
all the lower-dimensional cones that arise from the intersecting faces of
the subcones in an inclusion-exclusion formula.  
Later it was pointed out that it is possible to simplify the algorithm 
by computing with full-dimensional cones only, by making use of Brion's ``polarization trick'' \citep[see][Remark
4.3]{BarviPom}:
The computations with rational generating
functions are invariant with respect to the contribution of non-pointed cones
(cones containing a non-trivial linear subspace).  By operating in the dual space,
i.e., by computing with the polars of all cones, lower-dimensional cones can 
be safely discarded, because this is equivalent to discarding non-pointed
cones in the primal space.  The practical implementations also rely heavily on
this polarization trick.

\bigbreak

In practical implementations of Barvinok's algorithm, one observes that in the
hierarchy of cone decompositions, the index of the decomposed cones quickly
descends from large numbers to fairly low numbers. The
``last mile,'' i.e., decomposing many cones with fairly low index, creates
a huge number of unimodular cones and thus is the bottleneck of the whole
computation in many instances.  

The idea of this paper is to stop the decomposition when the index of a cone
is small enough, and to compute with generating functions for the integer
points in cones of small index rather than unimodular cones.  
When we try to implement this simple idea in Barvinok's algorithm, as outlined in
\autoref{dual}, we face a major difficulty, however:
Polarizing back a cone of small index
can create a cone of very large index, because determinants of $d\times d$
matrices are homogeneous of order~$d$.  

To address this difficulty, we avoid polarization altogether and perform the
signed decomposition in the primal space instead.  To avoid having to deal
with all the lower-dimensional subcones, we use the concept of
\emph{irrational decompositions} 
of rational polyhedra.  \cite{beck-sottile:irrational} introduced this notion
to give astonishingly simple proofs for three theorems of Stanley on
generating functions for the integer points in rational polyhedral cones.
Using the same technique, \cite{beck-haase-sottile:theorema} gave simplified
proofs of theorems of Brion and Lawrence\with Varchenko.  An irrational
decomposition of a polyhedron is a decomposition into polyhedra whose
proper faces do not contain any lattice points.  Counting formulas for
lattice points based on irrational decompositions therefore do not need to
take any inclusion-exclusion principle into account.

We give an explicit construction of a \emph{uniform irrational shifting
  vector} $\ve s$ for a cone $\ve v + K$ with apex $\ve v$ such that the
shifted cone $(\ve v + \ve s) + K$ has the same lattice points and contains no
lattice points on its proper faces (\autoref{construction}).  More strongly,
we prove that \emph{all 
  cones} appearing in the signed decompositions of $(\ve v + \ve s) + K$ in
Barvinok's algorithm contain no lattice points on their proper faces.
Therefore, discarding lower-dimensional cones is safe.  Despite its name, the
vector $\ve s$ only has \emph{rational coordinates}, so after shifting the cone by
$\ve s$, large parts of existing implementations of Barvinok's algorithm can
be reused to compute the irrational primal decompositions.

In \autoref{algorithm}, we show the precise algorithm.  
We also show that the same technique can be applied to the ``homogenized version'' of
Barvinok's algorithm that was proposed by \citet{latte2}. 

In \autoref{non-simplicial}, we extend the irrationalization technique to
non-simplicial cones.  This gives rise to an ``all-primal'' Barvinok
algorithm, where also triangulation of non-simplicial cones is performed in
the primal space.  This allows us to handle problems where the triangulation of
the dual cones is hard, e.g., in the case of cross polytopes. 

Finally, in \autoref{experiments}, we report on computational results.  Results on benchmark problems show that the new algorithms are
faster than the existing algorithms by orders of magnitude.  We also include 
results for problems that
could not previously be solved with Barvinok techniques.


\section{Barvinok's algorithm}
\renewcommand\textfraction{.5}

Let $P\subseteq\R^d$ be a rational polyhedron.  The \emph{generating function}
of~$P\cap\Z^d$ is defined as the formal Laurent series
\begin{displaymath}
  \tilde g_P(\ve z) = \sum_{\vealpha\in P\cap\Z^d} \ve z^{\vealpha}
\in \Z[[z_1,\dots,z_d, z_1^{-1},\dots,z_d^{-1}]],
\end{displaymath}
using 
the multi-exponent notation $\ve z^{\vealpha} = \prod_{i=1}^d z_i^{\alpha_i}$.
If $P$ is bounded, $\tilde g_P$ is a Laurent polynomial, which we
consider as a rational function~$g_P$.  If $P$ is not
bounded but is pointed (i.e., $P$ does not contain a straight line), there is
a non-empty open subset $U\subseteq\C^d$ such that the series converges
absolutely and uniformly on every compact subset of~$U$ to a rational
function~$g_P$.  If $P$ contains a straight line, we set $g_P \equiv 0$.
The rational function $g_P\in\Q(z_1,\dots,z_d)$ defined in this way is called
the \emph{rational generating function} of~$P\cap\Z^d$.

\begin{figure}[t]
  \begin{center} 
    \begin{minipage}[b]{.45\linewidth} 
      \ifpdf
    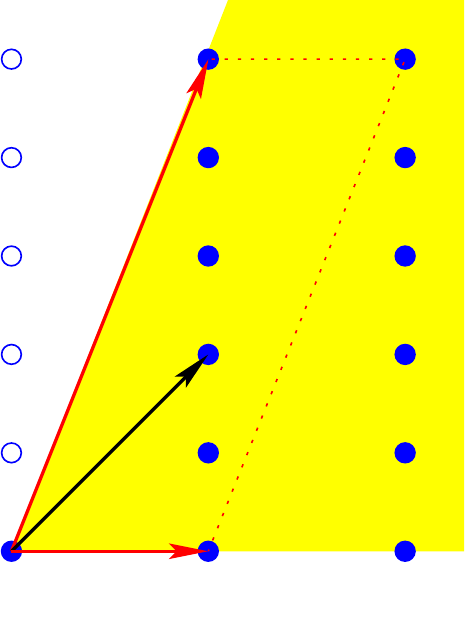
    \else
    \input{positive-decomp-1.pstex_t}
    \fi 
    \end{minipage}\qquad 
    \begin{minipage}[b]{.45\linewidth} 
      \ifpdf
    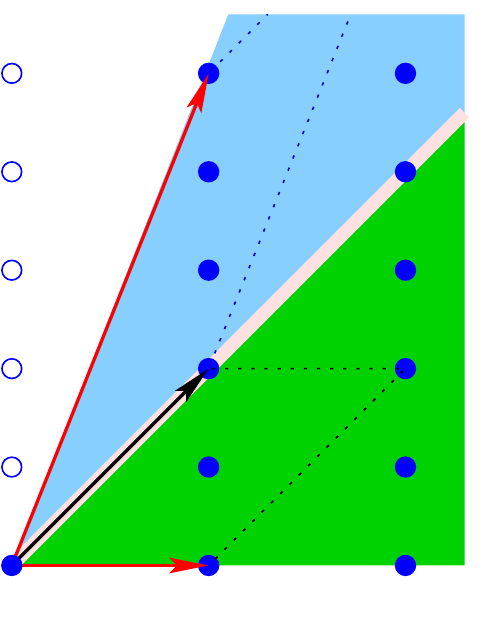
    \else
    \input{positive-decomp-2.pstex_t}
    \fi 
    \end{minipage} 
  \end{center} 
  \caption{A triangulation of the cone of index~$5$ generated by $\ve b^1$ and
    $\ve b^2$ into the two cones spanned by $\{\ve b^1,\ve w\}$ and $\{\ve
    b^2,\ve w\}$, having an index of $2$ and~$3$, respectively.  We have
    the inclusion-exclusion formula $g_{\cone\{\ve b_1,\ve b_2\}}(\ve z) = 
    g_{\cone\{\ve b_1,\ve w\}}(\ve z) + g_{\cone\{\ve b_2,\ve w\}}(\ve z) -
    g_{\cone\{\ve w\}}(\ve z)$; 
    here the one-dimensional cone spanned by $\ve w$ 
    needed to be subtracted.}
  \label{fig:positive-decomp}
\end{figure}

\begin{figure}[t]
  \begin{center} 
  \begin{minipage}[b]{.45\linewidth}
    \ifpdf
    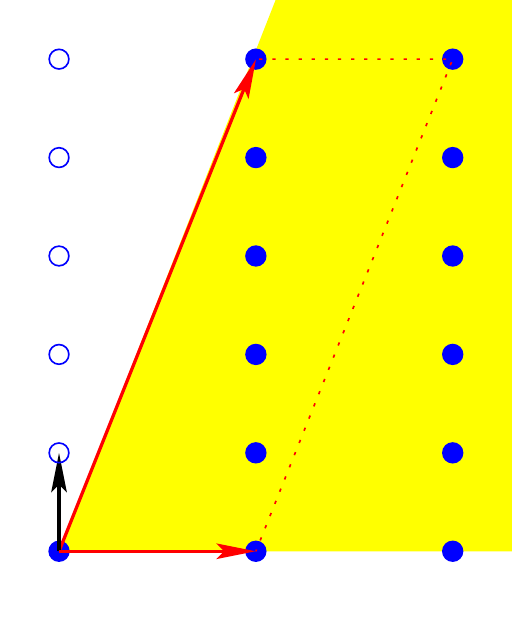
    \else
    \input{signed-decomp-1.pstex_t}
    \fi 
  \end{minipage}\qquad 
  \begin{minipage}[b]{.45\linewidth} 
    \ifpdf
    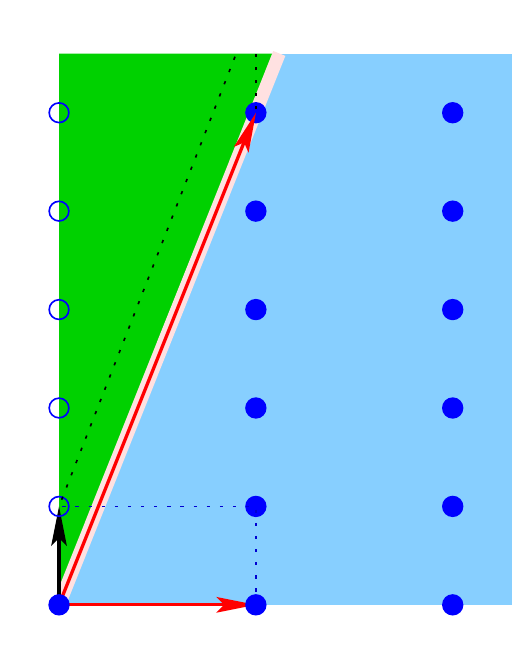
    \else
    \input{signed-decomp-3.pstex_t}
    \fi 
  \end{minipage} 
\end{center} 
  \caption{A signed decomposition of the cone of index~$5$ generated by $\ve b^1$ and
    $\ve b^2$ into the two unimodular cones spanned by $\{\ve b^1,\ve w\}$ and $\{\ve
    b^2,\ve w\}$.  We have the inclusion-\penalty0 exclusion formula $g_{\cone\{\ve b_1,\ve b_2\}}(\ve z)
    = g_{\cone\{\ve b_1,\ve w\}}(\ve z) - g_{\cone\{\ve b_2,\ve w\}}(\ve z) +
    g_{\cone\{\ve w\}}(\ve z)$.} 
  \label{fig:signed-decomp}
\end{figure}

Barvinok's algorithm computes the rational generating function of a
polyhedron~$P$.  It proceeds as follows.  By the theorem of \citet{Brion88},
the rational generating function of a 
polyhedron can be expressed as the sum of the rational generating functions of
the supporting cones of its vertices.  Let $\ve v_i\in\Q^d$ be a vertex of the
polyhedron~$P$.  Then the \emph{supporting cone} $\ve v_i + C_i$ of $\ve v_i$ is the
(shifted) polyhedral cone defined by $\ve v_i + \cone(P - \ve v_i)$.  Every supporting
cone $\ve v_i+C_i$ can be triangulated to obtain simplicial cones $\ve v_i+C_{ij}$.  
Let $K = \ve v+B\R_+^d$ be a simplicial full-dimensional cone, whose
\emph{basis vectors} $\ve b_1,\dots,\ve b_d$ (i.e., representatives
of its extreme rays) are
given by the columns of some matrix $B\in\Z^{d\times d}$.  We assume that the
basis vectors are primitive vectors of the standard lattice~$\Z^d$.
Then the \emph{index} of $K$ is defined to be $\ind K = \left|\det B\right|$; it
can also be interpreted as the cardinality of~$\Pi\cap\Z^d$, where $\Pi$ is
the \emph{fundamental parallelepiped} of $K$, i.e., the half-open parallelepiped
\begin{displaymath}
  \Pi = \ve v + \Bigl\{\, \textstyle\sum_{i=1}^d \lambda_i \ve b_i : 0 \leq \lambda_i <
    1\,\Bigr\}. 
\end{displaymath}
We remark that the set $\Pi\cap\Z^d$ can also be seen as a set of
representatives of the cosets of 
the lattice $B\Z^d$ in the standard lattice~$\Z^d$; we shall make use of this
interpretation in \autoref{dual}.  Barvinok's algorithm now
computes a \emph{signed decomposition} of the simplicial cone $K$ to
produce other simplicial cones with smaller index.  To this end, the algorithm
constructs a vector $\ve w\in\Z^d$ such that
\begin{equation}
  \label{eq:omega-vector-0}
  \ve w = \alpha_1 \ve b_1 + \dots + \alpha_d \ve b_d
  \quad\text{with}\quad \abs{\alpha_i} \leq \abs{\det B}^{-1/d} \leq 1.
\end{equation}
This can be accomplished using integer programming or lattice basis reduction.
The cone is then decomposed into cones spanned by $d$ vectors from the set
$\{\ve b_1,\dots,\ve b_d,\ve w\}$; each of the resulting cones then has an
index bounded above by $(\ind K)^{(d-1)/d}$.  In general, these cones form a
signed decomposition of~$K$ (see \autoref{fig:signed-decomp}); if $\ve w$ lies inside~$K$, 
they form a triangulation of~$K$ (see \autoref{fig:positive-decomp}).
The resulting cones and their intersecting proper faces (arising in an
inclusion-exclusion formula) are recursively processed, until
\emph{unimodular} cones, i.e., cones of index~$1$ are obtained.  
Finally, for a unimodular cone $\ve v+B\R_+^d$, the rational generating function can
be easily written down as 
\begin{equation}
  \frac{\ve z^{\ve a}}{\prod_{j=1}^d (1-\ve z^{\ve b_j})},
\end{equation}
where $\ve a$ is the unique integer point in the fundamental parallelepiped
of the cone.
We summarize Barvinok's algorithm below.\medbreak

\begin{algorithm}[Barvinok's original (primal) algorithm]~\smallskip\par
\label{algo:primal-barvi}
\noindent {\em Input:} 
A polyhedron $P\subset\R^d$ given by rational inequalities.

\noindent {\em Output:} 
The rational generating function for $P\cap\Z^d$ in the
form
\begin{equation}
  \label{eq:generating-function-0}
  g_P(\ve z) = \sum_{i\in I} \epsilon_i \frac{\ve
    z^{\ve a_i}}{\prod_{j=1}^d (1-\ve z^{\ve b_{ij}})}
\end{equation}
where $\epsilon_i\in\{\pm1\}$, $\ve a_i\in\Z^d$, and
$\ve b_{ij}\in\Z^d$.
\begin{enumerate}[\ 1.]
\item Compute all vertices $\ve v_i$ and corresponding 
  supporting cones $C_i$ of $P$.
\item Triangulate $C_i$ into simplicial cones $C_{ij}$, keeping track of all
  the intersecting proper faces.
\item Apply signed decomposition to the cones $\ve v_i + C_{ij}$ 
  to obtain unimodular cones $\ve v_i + C_{ijl}$, keeping track of all the
  intersecting proper faces. 
\item Compute the unique integer point~$\ve a_i$ in the fundamental parallelepiped
  of every resulting cone~$\ve v_i+C_{ijl}$.
\item Write down the formula~\eqref{eq:generating-function-0}.
\end{enumerate}
\end{algorithm}

The recursive decomposition of cones defines a \emph{decomposition tree}.  
Due to the descent of the indices in the signed decomposition
procedure, the following estimate holds for its depth:  
\begin{lemma}[\citealp{bar}]\label{lemma:barvi-depth}
  Let $B\R^d_+$ be a simplicial full-dimensional cone, whose basis is given by
  the columns of the matrix $B\in\Z^{d\times d}$.  Let $D=\abs{\det B}$.
  Then the depth of the decomposition tree is at most
  \begin{equation}\label{eq:barvi-depth}
    k(D) = \floor{1 + \frac{ \log_2 \log_2 D }{\log_2 \frac{d}{d-1} } }.
  \end{equation}
\end{lemma}
Because at each decomposition step at most $\mathrm O(2^d)$ cones are created
and the depth of the tree is doubly logarithmic in the index of the input
cone, Barvinok could obtain a polynomiality result \emph{in fixed dimension}: 
\begin{theorem}[\citealp{bar}]
  Let $d$ be fixed.  There exists a polynomial-time algorithm for computing 
  the rational generating function of a polyhedron $P\subseteq\R^d$ given by rational
  inequalities. 
\end{theorem}

Later the algorithm was improved by making use of Brion's ``polarization
trick'' \citep[see][Remark 4.3]{BarviPom}: The computations with rational
generating functions are invariant with respect to the contribution of
non-pointed cones (cones containing a non-trivial linear subspace). 
The reason is that the rational generating function of every non-pointed cone
is zero.  By operating in the dual space, i.e., by computing with the polars of all cones,
lower-dimensional cones can be safely discarded, because this is equivalent to
discarding non-pointed cones in the primal space.\medbreak  

\begin{algorithm}[Dual Barvinok algorithm]~\smallskip\par
\label{algo:dual-barvi}
\noindent {\em Input:} 
A polyhedron $P\subset\R^d$ given by rational inequalities.

\noindent {\em Output:} 
The rational generating function for $P\cap\Z^d$ in the
form
\begin{equation}
\label{eq:generating-function-1}
  g_P(\ve z) = \sum_{i\in I} \epsilon_i \frac{\ve
    z^{\ve a_i}}{\prod_{j=1}^d (1-\ve z^{\ve b_{ij}})}
\end{equation}
where $\epsilon_i\in\{\pm1\}$, $\ve a_i\in\Z^d$, and
$\ve b_{ij}\in\Z^d$.
\begin{enumerate}[\ 1.]
\item Compute all vertices $\ve v_i$ and corresponding 
  supporting cones $C_i$ of $P$.
\item Polarize the supporting cones $C_i$ to obtain $C_i^\circ$.
\item Triangulate $C_i^\circ$ into simplicial cones $C^\circ_{ij}$, discarding
  lower-dimensional cones. 
\item Apply Barvinok's signed decomposition to the cones $\ve v_i + C^\circ_{ij}$ 
  to obtain cones $\ve v_i + C^\circ_{ijl}$, stopping decomposition when a
  unimodular cone is obtained.  Discard all
  lower-dimensional cones. 
\item Polarize back $C^\circ_{ijl}$ to obtain cones $C_{ijl}$.
\item Compute the unique integer point~$\ve a_i$ in the fundamental parallelepiped
  of every resulting cone~$\ve v_i+C_{ijl}$.
\item Write down the formula~\eqref{eq:generating-function-1}.
\end{enumerate}
\end{algorithm}
 
This variant of the algorithm is much faster than the original algorithm
because in each step of the signed decomposition at most $d$, rather than
$\mathrm O(2^d)$, cones are created.  The practical
implementations LattE \citep{latte1} and \texttt{barvinok}
\citep{verdoolaege-et-al:counting-parametric} also rely heavily on this polarization trick.

\section{The Barvinok algorithm with stopped decomposition}
\label{dual}
We start out by introducing a first variant of Barvinok's algorithm that stops
decomposing cones before unimodular cones are reached.  
As we will see in the computational results
in \autoref{experiments}, already the simple modification that we propose can give a
significant improvement of the running time for some problems, at least in low
dimension. 
\medskip

\begin{algorithm}[Dual Barvinok algorithm with stopped decomposition]~\smallskip\par
\label{algo:dual-barvi}
\noindent {\em Input:} 
A polyhedron $P\subset\R^d$ given by rational inequalities; 
the maximum index $\ell$.

\noindent {\em Output:} 
The rational generating function for $P\cap\Z^d$ in the
form
\begin{equation}
  \label{eq:generating-function-d}
  g_P(\ve z) = \sum_{i\in I} \epsilon_i \frac{\sum_{\ve a\in A_i} \ve
    z^{\ve a}}{\prod_{j=1}^d (1-\ve z^{\ve b_{ij}})}
\end{equation}
where $\epsilon_i\in\{\pm1\}$, $A_i\subseteq\Z^d$ with $|A_i|\leq \ell$, and
$\ve b_{ij}\in\Z^d$.
\begin{enumerate}[\ 1.]
\item Compute all vertices $\ve v_i$ and corresponding 
  supporting cones $C_i$ of $P$.
\item Polarize the supporting cones $C_i$ to obtain $C_i^\circ$.
\item Triangulate $C_i^\circ$ into simplicial cones $C^\circ_{ij}$, discarding
  lower-dimensional cones. 
\item Apply Barvinok's signed decomposition to the cones $\ve v_i + C^\circ_{ij}$ 
  to obtain cones $\ve v_i + C^\circ_{ijl}$, stopping decomposition when a
  polarized-back cone 
  $C_{ijl} = (C^\circ_{ijl})^\circ$ has index at most~$\ell$.  Discard all
  lower-dimensional cones. 
\item Polarize back $C^\circ_{ijl}$ to obtain cones $C_{ijl}$.
\item Enumerate the integer points in the fundamental parallelepipeds 
  of all resulting cones~$\ve v_i+C_{ijl}$ to obtain the sets $A_i$.\label{algostep:enum}
\item Write down the formula~\eqref{eq:generating-function-d}.
\end{enumerate}
\end{algorithm}

As mentioned above, the integer points in the fundamental parallelepiped of a
cone $\ve v_i+B_{ijl}\R_+^d$ can be interpreted as representatives of the cosets of
the lattice $B_{ijl}\Z^d$ in the standard lattice $\Z^d$.  Hence they can be easily
enumerated in step~\ref{algostep:enum} by computing the Smith normal form of the generator matrix~$B_{ijl}$; see
Lemma~5.2 of \cite{barvinok-1993:exponential-sums}.  The Smith normal form
can be computed in polynomial time, even if the dimension is not fixed
\citep{kannan-bachem:79}.

We remark that both triangulation and signed decomposition are done in the dual
space, but the stopping criterion is the index of the polarized-back cones (in
the primal space).
The reason for this stopping criterion is that we wish to control the maximum
number of points in the fundamental parallelepipeds that need to be
enumerated.  Indeed, when the maximum $\ell$ is chosen as a constant or
polynomially in the input 
size, then \autoref{algo:dual-barvi} clearly runs in polynomial time (in fixed
dimension). 

Each step of Barvinok's signed decomposition reduces the index of the
decomposed cones.  When the index of a cone $C^\circ_{ijl}$ is~$\Delta$, in
the worst case the polarized-back cone $C_{ijl}$ has index
$\Delta^{d-1}$, where $d$ is the dimension.  If the dimension is too large,
the algorithm often needs to decompose cones down to a very low index or even index~1, so
the speed-up of the algorithm will be very limited.  This can be seen from the
computational results in \autoref{experiments}.

\section{Construction of a uniform irrational shifting vector}
\label{construction}

In this section, we will give an explicit construction of an \emph{irrational
  shifting vector} $\ve s$ for a simplicial cone $\ve v + K$ with 
apex $\ve v$ such that the shifted cone $(\ve v + \ve s) + K$ has the same
lattice points and contains no lattice points on its proper faces.  The
``irrationalization'' (or perturbation) will be \emph{uniform} in the sense that
also every cone arising during the Barvinok decomposition does not contain any
lattice points on its proper faces.  This will enable us to perform the Barvinok
decomposition in the primal space, discarding all lower-dimensional cones.

\begin{figure}[bt]
  \ifpdf
  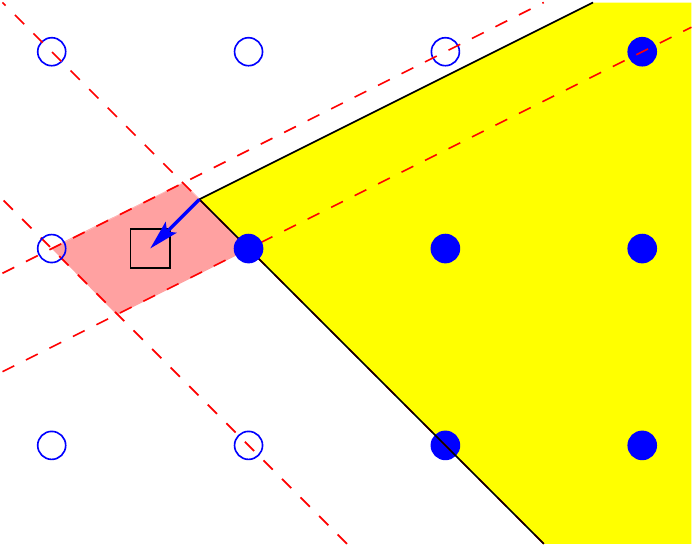
  \else
  \input{stability-region-3.pstex_t}
  \fi
  \caption{The stability region of a cone}
  \label{fig:stability}
\end{figure}

To accomplish this goal, we shall first describe a subset of the \emph{stability region}
of a cone $\ve v + K$ with apex 
at~$\ve v$, i.e., the set of apex points $\ve{\tilde v}$ such that $\ve{\tilde
  v}+K$ contains the same lattice points as $\ve v+K$; see Figure~\ref{fig:stability}.
\begin{lemma}[Stability cube]
  \label{lemma:stability}
  Let $\ve v + B\R^d_+$ be a simplicial full-dimensional cone with apex at $\ve
  v\in\Q^d$, whose basis is given by the columns of the matrix
  $B\in\Z^{d\times d}$.  Let $\ve b^*_1,\dots,\ve b^*_d$ be a basis of the dual
  cone, given by the columns of the matrix $B^* = -(B^{-1})^\top$.

  Let $D = \abs{\det B}$.    
  Let $\velambda\in\Q^d$ and $\vehatlambda\in\Q^d$ be defined by
  \begin{equation*}
    \label{eq:centerlambda}
    \lambda_i = \langle \ve b^*_i, \ve v\rangle
    \quad\text{and}\quad
    \hat\lambda_i = \frac1{D} \left(\floor{D
        \lambda_i} + \frac12\right) \quad\text{for $i=1,\dots,d$}.
  \end{equation*}
  Let 
  \begin{equation*}
    \label{eq:stability-length}
    \ve{\hat v} = -B\vehatlambda
    \quad\text{and}\quad \rho = \frac{1}{2D\cdot \max_{i=1}^d \norm{\ve b^*_i}_1}.
  \end{equation*}
  Then, for every $\ve{\tilde v}$ with $\norm{\ve{\tilde
      v} - \ve{\hat v} }_\infty < \rho$, the cone
  $\ve{\tilde v} + B\R^d_+$ contains the same integer points as the cone $\ve
    v+B\R^d_+$ and does not have integer points
    on its proper faces.
\end{lemma}
In the proof of the lemma, we will use of the inequality description
(H-representation) of the simplicial cone $\ve v + B\R^d_+$.  It is given by
the basis vectors of the dual cone:
\begin{equation}\label{eq:h-rep}
  \ve v + B\R^d_+ = \left\{\, \ve x\in\R^d : 
    \langle \ve b^*_i, \ve x \rangle \leq \langle \ve b^*_i, \ve v\rangle 
    \text{ for $i=1,\dots,d$} \,\right\}.
\end{equation}

\begin{proof}[Proof of \autoref{lemma:stability}]
  Let $\vetildelambda$ be defined by $\tilde\lambda_i = \langle \ve b^*_i,
  \ve{\tilde v}\rangle$.  Then we have
  \begin{equation}
    \label{eq:lambda-distance}
    \abs[big]{\tilde\lambda_i - \hat\lambda_i}
    \leq \norm{\ve b^*_i}_1 \cdot \norm{\ve{\tilde v} - \ve{\hat v}}_\infty
    < \norm{\ve b^*_i}_1 \cdot \rho 
    \leq \frac1{2D}.
  \end{equation}
  By~\eqref{eq:h-rep}, a point $\ve x\in\Z^d$ lies in the cone $\ve v + B\R^d_+$ if and
  only if 
  \begin{equation*}
    \label{eq:polarity-lambda}
    \langle \ve b^*_i, \ve x \rangle \leq \langle \ve b^*_i, \ve
    v\rangle = \lambda_i \quad\text{for $i=1,\dots,d$.}
  \end{equation*}
  Likewise, $\ve x\in\ve{\tilde v} + B\R^d_+$ if and only if
  \begin{equation*}
    \label{eq:polarity-lambda-tilde}
    \langle \ve b^*_i, \ve x \rangle \leq \langle \ve b^*_i, \ve{\tilde v}
    \rangle = \tilde\lambda_i \quad\text{for $i=1,\dots,d$.}
  \end{equation*}
  Note that for $\ve x\in\Z^d$, the left-hand sides of both inequalities are an
  integer multiple of $\frac1D$. Therefore, we obtain equivalent statements by
  rounding down the right-hand sides to integer multiples of $\frac1D$.
  For the right-hand side of~\eqref{eq:polarity-lambda-tilde} we have by~\eqref{eq:lambda-distance}
  \begin{subequations}
    \label{eq:included-in-open-stability}
    \begin{align}
      \tilde\lambda_i &= \hat\lambda_i + (\tilde \lambda_i - \hat\lambda_i) <
      \frac 1D \left( \floor{D  \lambda_i} + \frac12 \right)
      + \frac1{2D} = \frac 1D \floor{D \lambda_i} + 1, \\
      \tilde\lambda_i &= \hat\lambda_i + (\tilde \lambda_i - \hat\lambda_i) >
      \frac 1D \left( \floor{D \lambda_i} + \frac12 \right) - \frac1{2D}
      = \frac 1D \floor{D \lambda_i},
    \end{align}
  \end{subequations}
  so $\lambda_i$ and $\tilde\lambda_i$ are rounded down to the same
  value $\frac1D\floor{D\lambda_i}$.  Thus, the cone $\ve{\tilde v} + B\R^d_+$
  contains the same integer points as the cone $\ve v+B\R^d_+$.  Moreover,
  since the inequalities~\eqref{eq:included-in-open-stability} are strict, the
  cone $\ve{\tilde v} + B\R^d_+$ does not have integer points on its proper
  faces.
\end{proof}
For non-simplicial cones, we will give an algorithmic construction for a
stability cube in \autoref{non-simplicial}.\smallbreak

Next we make use of the estimate for the depth of the decomposition tree in Barvinok's
algorithm given in \autoref{lemma:barvi-depth}.  
On each level of the decomposition, the entries in the basis matrices can
grow, but not by much.  We obtain:
\begin{lemma}\label{lemma:basis-bound}
  Let $B\R^d_+$ be a simplicial full-dimensional cone, whose basis is given by
  the columns of the matrix $B\in\Z^{d\times d}$.  Let $D=\abs{\det B}$.  Let
  $C\in\Z_+$ be a number  such that $\abs{B_{i,j}} \leq C$. 

  Then all the basis matrices $\bar B$ of the cones that appear in the
  recursive signed decomposition procedure of Barvinok's
  algorithm applied to $B\R^d_+$ have entries bounded above by $d^{k(D)} C$,
  where $k(D)$ is defined by~\eqref{eq:barvi-depth}.
\end{lemma}
\begin{proof}
  Given a cone spanned by the columns~$\ve b_1,\dots,\ve b_d$ of the
  matrix~$B$, Barvinok's algorithm constructs a vector $\ve w\in\Z^d$ such
  that
  \begin{equation}
    \label{eq:omega-vector}
    \ve w = \alpha_1 \ve b_1 + \dots + \alpha_d \ve b_d
    \quad\text{with}\quad \abs{\alpha_i} \leq \abs{\det B}^{-1/d} \leq 1.
  \end{equation}
  Thus ${\norm{\ve w}}_\infty \leq d C$.  
  The cone is then decomposed into cones spanned by $d$ vectors from the set
  $\{\ve b_1,\dots,\ve b_d,\ve w\}$.  Thus the entries in the corresponding
  basis matrices are bounded by $dC$.  The result follows then by
  \autoref{lemma:barvi-depth}. 
\end{proof}

If we can bound the entries of an integer matrix with non-zero determinant,
we can also bound the entries of its inverse.
\begin{lemma}\label{lemma:inverse-bound}
  Let $B\in\Z^{d\times d}$ be a matrix with $\abs{B_{i,j}} \leq C$.  Let
  $D=\abs{\det B}$. Then the
  absolute values of the entries of $B^{-1}$ are bounded above by 
  \begin{displaymath}
    \frac1D
    (d-1)! C^{d-1}.
  \end{displaymath}
\end{lemma}
\begin{proof}
  We have $\abs{(B^{-1})_{k,l}} = \frac1D \abs{\det B_{(k,l)}}$, where $B_{(k,l)}$ is
  the matrix obtained from deleting the $k$-th row and $l$-th column from~$B$. 
  Now the desired estimate follows from a formula for $\det B_{(k,l)}$ and from
  $\abs{B_{i,j}} \leq C$. 
\end{proof}

Thus, we obtain a bound on the norm of the basis vectors
of the polars of all cones occurring in the signed decomposition procedure of
Barvinok's algorithm. 
\begin{corollary}[A bound on the dual basis vectors]
  \label{lemma:dual-basis-bound}
  Let $B\R^d_+$ be a simplicial full-dimensional cone, whose basis is given by
  the columns of the matrix $B\in\Z^{d\times d}$.  Let $D=\abs{\det B}$.  Let
  $C$ be a number such that $\abs{B_{i,j}} \leq C$. 

  Let $\bar B^* = -(\bar B^{-1})^\top$ be the basis matrix of the polar
  of an arbitrary cone $\bar B\R_+^d$ that appears in the recursive signed decomposition
  procedure applied to $B\R^d_+$.
  Then, for every column vector $\ve{\bar b}^*_i$ of $\bar B^*$ we have the estimate
  \begin{equation}
    \label{eq:dual-basis-bound}
    \norm{\det{\bar B} \cdot \ve b^*_i}_\infty
    \leq (d-1)! \paren{ d^{k(D)} C }^{d-1} =: L
  \end{equation}
  where $k(D)$ is defined by~\eqref{eq:barvi-depth}.
\end{corollary}
\begin{proof}
  By \autoref{lemma:basis-bound}, the entries of $\bar B$ are bounded above by
  $d^{k(D)} C$.  Then the result follows from \autoref{lemma:inverse-bound}.
\end{proof}

The construction of the ``irrational'' shifting vector is based on the
following lemma.

\begin{lemma}[The irrational lemma]
  \label{lemma:irrational}
  Let $M\in\Z_+$ be an integer.  Let
  \begin{equation}
    \label{eq:irrational-vector}
    \ve q = \left( \frac1{2M}, \frac1{(2M)^2}, \dots, \frac1{(2M)^d} \right).
  \end{equation}
  Then $\langle \ve c,\ve q\rangle \notin\Z$ for every
  $\ve c\in\Z^d\setminus\{\ve0\}$ with $\norm{\ve c}_\infty < M$.
\end{lemma}
\begin{proof}
  Follows from the principle of representations of rational numbers in a
  positional system of base $2M$. 
\end{proof}

\begin{theorem}
  \label{th:uniform-irr}
  Let $\ve v + B\R^d_+$ be a simplicial full-dimensional cone with apex at $\ve
  v\in\Q^d$, whose basis is given by the columns of the matrix
  $B\in\Z^{d\times d}$.  Let $D=\abs{\det B}$, $C$ be a number such that
  $\abs{B_{i,j}} \leq C$ and let $\ve{\hat v}\in\Q^d$ and $\rho\in\Q_+$ be the
  data from \autoref{lemma:stability} describing the stability cube of~$\ve
  v + B\R^d_+$.  Let $0<r\in\Z$ such that $r^{-1} < \rho$.
  Using 
  \begin{equation}
    k = \floor{1 + \frac{ \log_2 \log_2 D }{\log_2 \frac{d}{d-1} } },
  \end{equation}
  $L = (d-1)! ( d^{k} C )^{d-1} $, and $M = 2L$, 
  define
  \begin{displaymath}
    \ve s = \frac1r \cdot\left(\frac1{(2M)^1},\frac1{(2M)^2}, \dots, \frac1{(2M)^d} \right).
  \end{displaymath}
  Finally let $\ve{\tilde v} = \ve{\hat v} + \ve s$.
  \begin{enumerate}[\rm (i)]
  \item We have $(\ve{\tilde v} + B\R^d_+)\cap\Z^d = (\ve{v} + B\R^d_+)\cap\Z^d$,
    i.e., the shifted cone has the same set of integer
    points as the original cone.
  \item The shifted cone $\ve{\tilde v} + B\R^d_+$ contains no lattice points on
    its proper faces. 
  \item More strongly, \emph{all cones} appearing in the signed decompositions
    of the shifted cone $\ve{\tilde v} + B\R^d_+$ in Barvinok's algorithm contain no lattice
    points on their proper faces.
  \end{enumerate}
\end{theorem}

\begin{proof}
  \emph{Part (i).}  This follows from \autoref{lemma:stability}
  because $\ve{\tilde v}$ clearly lies in the open stability cube.
  \smallbreak
  
  \emph{Parts (ii) and (iii).}
  Every cone appearing in the course of Barvinok's signed decomposition
  algorithm has the same apex $\ve{\tilde v}$ as the input cone and a basis
  $\bar B\in\Z^{d\times d}$ with $\abs{\det \bar B} \leq D$.
  Let such a $\bar B$ be fixed and denote by $\ve{\bar b}^*_i$ the columns of
  the dual basis matrix $\bar B^* = -(\bar B^{-1})^\top$.  
  Let $\ve z\in\Z^d$ be an arbitrary integer
  point.  We shall show that $\ve z$ is not on any of the facets of the cone,
  i.e., 
  \begin{equation}\label{eq:not-on-facet}
    \langle \ve{\bar b}^*_i, \ve z-\ve{\tilde v}\rangle\neq0\qquad\text{for $i=1,\dots,d$.}
  \end{equation}
  Let $i\in\{1,\dots,d\}$ arbitrary.  We will show~\eqref{eq:not-on-facet} by
  proving that
  \begin{equation}
    \label{eq:not-on-facet-2}
    \langle \det{\bar B}\cdot \ve{\bar b}^*_i, \ve{\tilde v}\rangle \notin\Z.
  \end{equation}
  Clearly, if \eqref{eq:not-on-facet-2} holds, we have $\langle \ve{\bar
    b}^*_i, \ve{\tilde v}\rangle \notin(\det \bar B)^{-1} \Z$.  But since $\langle
  \ve{\bar b}^*_i, \ve z\rangle \in (\det \bar B)^{-1} \Z$, we have 
  $\langle \ve{\bar b}^*_i, \ve z-\ve{\tilde v}\rangle\notin\Z$; in particular it is
  nonzero, which proves~\eqref{eq:not-on-facet}.

  To prove~\eqref{eq:not-on-facet-2}, let $\ve c=\det\bar B\cdot \ve{\bar
    b}^*_i$.  By \autoref{lemma:dual-basis-bound}, we have 
  $\norm{\ve c}_\infty \leq L < M$.
  Now \autoref{lemma:irrational} 
  gives $\langle\ve c, \ve s\rangle \notin \frac1r \Z$.
  Observing that by the definition~\eqref{eq:centerlambda}, we have
  \begin{displaymath}
    \langle \ve c, \ve{\hat v} \rangle 
    = \langle \ve c, -B\vehatlambda \rangle \in \tfrac1r \Z.
  \end{displaymath}
  Therefore, we have $\langle \ve c, \ve{\tilde v}\rangle
  = \langle \ve c, \ve{\hat v} + \ve s \rangle \notin\frac1r\Z$.
  This proves~\eqref{eq:not-on-facet-2}, and thus completes the proof.
\end{proof}

\section{The irrational algorithms}
\label{algorithm}

The following is our variant of the Barvinok algorithm.\medskip

\begin{algorithm}[Primal irrational Barvinok algorithm]~\smallskip\par
\label{algo:primal-irr-barvi}
\noindent {\em Input:} 
A polyhedron $P\subset\R^d$ given by rational inequalities; 
the maximum index $\ell$.

\noindent {\em Output:} 
The rational generating function for $P\cap\Z^d$ in the
form
\begin{equation}
  \label{eq:generating-function}
  g_P(\ve z) = \sum_{i\in I} \epsilon_i \frac{\sum_{\ve a\in A_i} \ve
    z^{\ve a}}{\prod_{j=1}^d (1-\ve z^{\ve b_{ij}})}
\end{equation}
where $\epsilon_i\in\{\pm1\}$, $A_i\subseteq\Z^d$ with $|A_i|\leq \ell$, and
$\ve b_{ij}\in\Z^d$.
\begin{enumerate}[\ 1.]
\item Compute all vertices $\ve v_i$ and corresponding 
    supporting cones $C_i$ of $P$.
  \item Polarize the supporting cones $C_i$ to obtain $C_i^\circ$.
  \item Triangulate $C_i^\circ$ into simplicial cones $C^\circ_{ij}$, discarding
    lower-dimensional cones. 
  \item Polarize back $C^\circ_{ij}$ to obtain simplicial cones $C_{ij}$.
  \item Irrationalize all cones by computing new apex vectors $\ve{\tilde
      v}_{ij}\in\Q^d$ from $\ve v_i$ and $C_{ij}$ as in
    \autoref{th:uniform-irr}.\label{algostep:irr} 
  \item Apply Barvinok's signed decomposition to the cones $\ve{\tilde
      v}_{ij} + C_{ij}$, discarding lower-dimensional cones, until all cones
    have index at most $\ell$. 
  \item Enumerate the integer points in the fundamental parallelepipeds 
    of all resulting cones to obtain the sets $A_i$.
  \item Write down the formula~\eqref{eq:generating-function}.
  \end{enumerate}
\end{algorithm}
\begin{theorem}
  \autoref{algo:primal-irr-barvi} is correct and runs in polynomial time when
  the dimension $d$ is fixed and the maximum index $\ell$ is bounded by a
  polynomial in the input size.
\end{theorem}
\begin{proof}
  This is an immediate consequence of the analysis of Barvinok's algorithm.
  The irrationalization (step \ref{algostep:irr} of the algorithm) increases
  the encoding length of the apex vector only by a polynomial amount, because
  the dimension $d$ is fixed and the depth $k$ only depends doubly logarithmic
  on the initial index of the cone.
\end{proof}

The same technique can also be applied to the ``homogenized version'' of
Barvinok's algorithm that was proposed by \citet{latte2}; see also
\citet[Algorithm~11]{latte1}. 

\medskip
  
\begin{algorithm}[Irrational homogenized Barvinok algorithm]~\smallskip\par
\label{algo:homog-irr-barvi}
\noindent {\em Input:} 
A polyhedron $P\subset\R^d$ given by rational inequalities in the form $A\ve
x\leq\ve b$; the maximum index $\ell$.

\noindent {\em Output:} 
A rational generating function in the form~\eqref{eq:generating-function} for
the integer points in the \emph{homogenization} of~$P$, i.e., the cone
\begin{equation}
  \label{eq:homogenization}
  C=\set{\, (\xi\ve x, \xi): \ve x\in P,\, \xi\in\R_+\,}.
\end{equation}
\begin{enumerate}[\ 1.]
\item Consider the inequality description for~$C$; it is given by $A\ve x -
  \ve b\xi \leq 0$.  The polar $C^\circ$ then has the rays $(A_{i,\cdot},
  -b_i)$, $i=1,\dots,m$.
\item Triangulate $C^\circ$ into simplicial cones $C^\circ_{j}$, discarding
  lower-dimensional cones. 
\item Polarize back the cones $C^\circ_j$ to obtain simplicial cones $C_j$. 
\item Irrationalize the cones $C_j$ to obtain shifted cones $\ve{\tilde v}_j+C_j$.
\item Apply Barvinok's signed decomposition to the cones $\ve{\tilde
    v}_j + C_j$, discarding lower-dimensional cones, until all cones have
  index at most $\ell$. 
\item Write down the generating function.
\end{enumerate}
  
\end{algorithm}

\section{Extension to the non-simplicial case}
\label{non-simplicial}
For polyhedral cones with few rays and many facets, it is usually much faster
to perform triangulation in the primal space than in the dual space,
cf.~\cite{bueler-enge-fukuda-2000:exact-volume}.  In this section, we show 
how to perform both Barvinok decomposition and triangulation in the primal
space.

The key idea is to use linear programming to compute a subset of the
stability region of the non-simplicial cones.
\begin{lemma}
  \label{lemma:general-stability-region}
  There is a polynomial-time algorithm that, given the vertex $\ve v\in\Q^d$
  and the facet vectors $\ve b^*_i\in\Z^d$, $i=1,\dots,m$, of a
  full-dimensional polyhedral cone~$C=\ve v+B\R^n_+$, where $n\geq d$,
  computes a point $\ve{\hat v}\in\Q^d$ and a positive scalar $\rho\in\Q$
  such that for every $\ve{\tilde v}$ in the open cube with $\norm{\ve{\tilde
      v} - \ve{\hat v}}_{\infty} < \rho$, the cone $\ve{\tilde v}+B\R^n_+$
  has no integer points on its proper faces and contains the same integer
  points as $\ve v+B\R^n_+$. 
\end{lemma} 
\begin{proof}
  We maximize $\rho$ subject to the linear inequalities
  \begin{subequations} 
    \begin{align} 
      \langle \ve b^*_i, \ve{\hat v}\rangle + \norm{\ve b^*_i}_1 \rho
      &\leq 
      \floor{\langle \ve b^*_i, \ve v\rangle} + 1,\label{eq:stablp-a}\\
      - \langle \ve b^*_i, \ve{\hat v}\rangle + \norm{\ve b^*_i}_1 \rho
      &\leq - \floor{\langle \ve b^*_i, \ve v\rangle},\label{eq:stablp-b} 
    \end{align} 
  \end{subequations} 
  where $\ve{\hat v}\in\R^d$ and $\rho\in\R_+$.  We can solve this linear
  optimization problem in polynomial time.  Let $(\ve{\hat v},\rho)$ be an
  optimal solution.  Let $\ve{\tilde v}\in\R^d$ with $\norm{\ve{\tilde
      v} - \ve{\hat v}}_{\infty} < \rho$.
  Let $\ve x\in(\ve{\tilde v} + B\R^d_+)\cap\Z^d$.  Then we have for every $i\in\{1,\dots,m\}$
  \begin{align*}
    \langle \ve b^*_i, \ve x \rangle & \leq \langle \ve b^*_i, \ve{\tilde v}\rangle\\
    &= \langle \ve b^*_i, \ve{\hat v} \rangle
    + \langle \ve b^*_i, \ve{\tilde v} - \ve{\hat v} \rangle\\
    &\leq \langle \ve b^*_i, \ve{\hat v} \rangle
    + \norm{\ve b^*_i}_1 \norm{\ve{\tilde v} - \ve{\hat v}}_\infty \\
    & < \langle \ve b^*_i, \ve{\hat v} \rangle + \norm{\ve b^*_i}_1 \rho \\
    & \leq \floor{\langle \ve b^*_i, \ve v\rangle} + 1
    && \text{by \eqref{eq:stablp-a}.}
  \end{align*}
  Because $\langle \ve b^*_i, \ve x \rangle$ is integer, we actually have
  \begin{math}
    \langle \ve b^*_i, \ve x \rangle \leq \floor{\langle \ve b^*_i, \ve v\rangle}
  \end{math}.
  Thus, $\ve x$ lies in the cone $\ve v + B\R^d_+$.
  Conversely, let $\ve x\in(\ve{v} + B\R^d_+)\cap\Z^d$.
  Then, for  every $i\in\{1,\dots,m\}$, we have
  \begin{math}
    \langle \ve b^*_i, \ve x \rangle \leq \langle \ve b^*_i, \ve v\rangle
  \end{math}.
  Since $\ve x\in\Z^d$, we can round down the right-hand side and obtain
   \begin{align*}
    \langle \ve b^*_i, \ve x \rangle & \leq \floor{\langle \ve b^*_i, \ve
      v\rangle} \\
    &\leq \langle\ve b^*_i, \ve{\hat v}\rangle - \norm{\ve b^*_i}_1 \rho &&
    \text{by \eqref{eq:stablp-b}}\\
    &<  \langle\ve b^*_i, \ve{\hat v}\rangle - \norm{\ve b^*_i}_1
    \norm{\ve{\tilde v} - \ve{\hat v}}_\infty \\
    &\leq \langle\ve b^*_i, \ve{\hat v}\rangle + \langle \ve b^*_i, \ve{\tilde
      v} - \ve{\hat v}\rangle \\
    &= \langle\ve b^*_i, \ve{\tilde v} \rangle.
  \end{align*}
  Thus, $\ve x\in\ve{\tilde v} + B\R^d_+$.  Moreover, since the inequality
  is strict, $\ve x$ does not lie on the face $\langle \ve b^*_i, \ve x
  \rangle = \langle \ve b^*_i, \ve{\tilde v}\rangle $ of the cone $\ve{\tilde v} +
  B\R^d_+$. 
\end{proof}~

\begin{lemma}[Bound on the index of all subcones]
  Let $\ve b_i\in\Z^d$, $i=1,\dots,n$, be the generators of a full-dimensional
  polyhedral cone~$K\subseteq\R^d$.  Then the cones of any
  triangulation of~$K$ have an index bounded by
  \begin{equation}
    \label{eq:triangulated-index-bound}
    D = \left( \max\nolimits_{i=1}^n \norm{\ve b_i}^2 \right)^{n/2}.
  \end{equation}
\end{lemma}
\begin{proof}
  Let $B\in\Z^{d\times d}$ be the generator matrix of a full-dimensional cone
  of a triangulation of~$K$; then the columns $B$ form a subset
  $\{\ve b_{i_1},\dots,\ve b_{i_d}\} \subseteq \{\ve b_1,\dots,\ve b_n\}$.  
  Therefore
  \begin{displaymath}
    \abs{\det B} \leq \prod_{k=1}^d \norm{\ve b_{i_k}} \leq \left( \max\nolimits_{i=1}^n
      \norm{\ve b_i}^2 \right)^{n/2},
  \end{displaymath}
  giving the desired bound.
\end{proof}

With these preparations, the following corollary is immediate.

\begin{corollary}
  \label{cor:uniform-irr-general}
  Let $\ve v + B\R_+^n$ be a full-dimensional polyhedral cone with apex at~$\ve
  v\in\Q^d$,  whose basis is given by the columns of the matrix
  $B\in\Z^{d\times d}$.  Let $D$ be defined
  by~\eqref{eq:triangulated-index-bound}.  Let $\ve{\hat v}\in\Q^d$ and
  $\rho\in\Q_+$ be the data from \autoref{lemma:general-stability-region}
  describing the stability cube of $\ve v + B\R_+^n$.  Using these data,
  construct $\ve{\tilde v}$ as in \autoref{th:uniform-irr}.
  Then the assertions of \autoref{th:uniform-irr} hold.  
\end{corollary}
~

\begin{algorithm}[All-primal irrational Barvinok algorithm]~\smallskip\par
  \label{algo:all-primal-irr-barvi}
  \noindent {\em Input:} 
  A polyhedron $P\subset\R^d$ given by rational inequalities; 
  the maximum index $\ell$.\par
  \noindent {\em Output:} 
  The rational generating function for $P\cap\Z^d$ in the
  form~\eqref{eq:generating-function}.
  \begin{enumerate}[\ 1.]
  \item Compute all vertices $\ve v_i$ and corresponding 
    supporting cones $C_i$ of $P$.
  \item Irrationalize all cones by computing new apex vectors $\ve{\tilde
      v}_{i}\in\Q^d$ from $\ve v_i$ by
    \autoref{cor:uniform-irr-general}.
  \item Triangulate $\ve{\tilde v}_i + C_i$ into simplicial cones $\ve{\tilde
      v}_i + C_{ij}$, discarding lower-dimensional cones. 
  \item Apply Barvinok's signed decomposition to the cones $\ve{\tilde
      v}_{i} + C_{ij}$, until all cones have index at most $\ell$.
  \item Enumerate the integer points in the fundamental parallelepipeds 
    of all resulting cones to obtain the sets $A_i$.
  \item Write down the formula~\eqref{eq:generating-function}.
  \end{enumerate}
\end{algorithm}

\section{Computational experiments}
\label{experiments}
\renewcommand\textfraction{0}

Algorithms~\ref{algo:primal-irr-barvi} and~\ref{algo:all-primal-irr-barvi}
have been implemented in a new version of the software package LattE, derived
from the official LattE release~1.2 \citep{latte-1.2}.  The new version,
called LattE macchiato, is
freely available on the Internet \citep{latte-macchiato}.  In this section, we
discuss some implementation details and show the results of first
computational experiments. 

\subsection{Two substitution methods}

When the generating function $g_P$ has been computed, the number of lattice
points can be obtained by evaluating $g_P(\ve 1)$.  However, $\ve 1$ is a
pole of every summand of the expression
\begin{displaymath}
  g_P(\ve z) = \sum_{i\in I} \epsilon_i \frac{\sum_{\ve a\in A_i} \ve
    z^{\ve a}}{\prod_{j=1}^d (1-\ve z^{\ve b_{ij}})}.
\end{displaymath}

The method implemented in LattE 1.2 (see \citealp{latte1}) is to use the
\emph{polynomial substitution}
\begin{displaymath}
  \ve z = ((1+s)^{\lambda_1}, \dots, (1+s)^{\lambda_d}).
\end{displaymath}
for a suitable vector $\velambda$.  Then the constant coefficient of the
Laurent expansion of every summand about $s=0$ is computed using polynomial
division.  The sum of all the constant coefficients finally gives the number
of lattice points.

Another method from the literature (see, for instance, \citealp{BarviPom}) is to 
use the \emph{exponential substitution}
\begin{displaymath}
  \ve z = (\exp\{\tau \lambda_1\}, \dots, \exp\{\tau \lambda_d\})
\end{displaymath}
for a suitable vector $\velambda$.  By letting $\tau\to0$, one then obtains
the formula 
\begin{equation}
  \label{eq:count-with-todd}
  |P\cap\Z^d| = \sum_{i\in I} \epsilon_i 
  \sum_{k=0}^d  \frac{\td_{d-k}( \langle \velambda,
    \ve b_{i1} \rangle, \dots, \langle \velambda, \ve b_{id}\rangle)}
  {k!\cdot \langle \velambda,
    \ve b_{i1} \rangle \cdots \langle \velambda, \ve b_{id}\rangle}
  \sum_{\ve a\in A_i} {\langle \velambda, \ve a\rangle}^k,
\end{equation}
where $\td_{d-k}$ is the so-called Todd polynomial.
In LattE macchiato, the exponential substitution method has been
implemented in addition to the existing polynomial substitution;
see \citet{latte-flavors} for implementation details.

\subsection{Implementation details}

We enumerate the lattice points in the fundamental parallelepiped by computing
the Smith normal form of the generator matrix~$B$; see Lemma~5.2 of
\cite{barvinok-1993:exponential-sums}.\footnote{The author wishes 
  to thank Susan Margulies for prototyping the enumeration code.}
For computing Smith normal
forms, we use the library \citet{lidia-2.2.0}.  For solving the linear program 
in \autoref{lemma:general-stability-region}, we use the implementation of the
revised dual simplex method in exact rational arithmetic in \texttt{cddlib},
version~0.94a \citep{cddlib-094a}. 
All other computations are
done using the libraries NTL, version~5.4 \citep{ntl-5.4} and
\citet{gmp-4.1.4} for providing exact integer and rational arithmetic.

\subsection{Evaluation of variants of the algorithms}
We compare the variants of the algorithms using test instances that can also
be solved without the proposed irrationalization techniques.  We consider
the test instances \texttt{hickerson-12}, \texttt{hickerson-13}, and
\texttt{hickerson-14}, related to 
the manuscript by \citet{hickerson-1991}.  They
describe simplices in $\R^6$ and $\R^7$ that contain $38$, $14$, and $32$ integer points,
respectively.  The examples are good test cases for our algorithms because 
the vertices and cones are trivially computed, and all computation time is
spent in the Barvinok decomposition.  
We show the results in \autoref{hickerson-12},
\autoref{hickerson-13}, and \autoref{hickerson-14}. 
The tables show results for the following methods:
\begin{enumerate}[\ 1.]
\item Methods without irrationalization, using polarization to avoid
  computing with lower-dimensional cones: 
  \begin{enumerate}[(a)]
  \item LattE 1.2 \citep{latte-1.2}, decomposing down to unimodular cones in the dual space
    (\autoref{algo:dual-barvi} with $\ell=1$).
  \item Likewise, but using the implementation in the
    library \texttt{barvinok} by \citet{barvinok-0.21}, version~0.21.
  \item LattE macchiato, decomposing cones in the dual space, until all
    cones in the primal space have at most index $\ell$ (\autoref{algo:dual-barvi}), 
    then using polynomial substitution.  We show the results for different
    values of~$\ell$.
  \item Likewise, but using exponential substitution.
  \end{enumerate}
\item Methods with irrationalization, performing triangulation in the dual space
  and Barvinok decomposition in the primal space (\autoref{algo:primal-irr-barvi}):
  \begin{enumerate}[(a)]
  \item LattE macchiato with polynomial substitution.
  \item LattE macchiato with exponential substitution.
  \end{enumerate}
\end{enumerate}
The table shows computation times in CPU seconds on a PC with a Pentium~M
processor with 1.4\,GHz.  It also shows the total number of simplicial cones
created in the decomposition, using the different variants of LattE; note that
we did not measure the number of simplicial cones that the library
\texttt{barvinok} produced.

\begin{table}[tp]
  \small
  \caption{Results for \texttt{hickerson-12}}
  \label{hickerson-12}
  \begin{center}
    \def~{\hphantom{0}}
    \begin{tabular}{cc*4{>{\sffamily}c}c*2{>{\sffamily}c}}
      \toprule
      & \multicolumn{5}{c}{Without irrationalization} 
      & \multicolumn{3}{c}{With irrationalization} \\
      \cmidrule{2-6} \cmidrule(l){7-9}
      & & \multicolumn{4}{c}{Time (s)} 
      & & \multicolumn{2}{c}{Time (s)} \\
      \cmidrule{3-6} \cmidrule{8-9}
      & & & & \multicolumn{2}{c}{LattE macchiato}
      & & \multicolumn{2}{c}{LattE macchiato} \\
      \cmidrule{5-6} \cmidrule{8-9}
      \smash[t]{\begin{tabular}[b]{c}Max.\\index\end{tabular}} 
      & Cones & \multicolumn{1}{c}{\smash[t]{\begin{tabular}[b]{@{}c@{}}LattE\\
            v\,1.2\end{tabular}} }
       & \multicolumn{1}{c}{\smash[t]{\begin{tabular}[b]{@{}c@{}}\texttt{barv.}\\
             v\,0.21\end{tabular}} }
      & \multicolumn{1}{c}{Poly} & \multicolumn{1}{c}{Exp} & Cones &
      \multicolumn{1}{c}{Poly} & \multicolumn{1}{c}{Exp} \\ 
      \midrule
      ~~~1 & 11625 & 17.9 & 11.9            & 10.0          & 16.7          & 7929         & 7.8          & 12.7          \\
      ~~10 & ~4251 & {}   & {}              & ~6.9          & ~7.0          & ~803         & 1.9          & ~1.6          \\
      ~100 & ~~980 & {}   & {}              & ~6.9          & ~2.1          & ~~84         & \textbf{1.3} & ~0.3 \\
      ~200 & ~~550 & {}   & {}              & ~9.1          & ~1.5          & ~~76         & 1.3          & ~0.3          \\
      ~300 & ~~474 & {}   & {}              & ~9.9          & ~1.4          & ~~58         & 1.4          & ~0.3          \\
      ~500 & ~~410 & {}   & {}              & 11.7          & ~1.3          & ~~42         & 1.6          & ~0.3          \\
      1000 & ~~130 & {}   & {}              & ~7.2          & ~0.7          & ~~22         & 1.7          & ~\textbf{0.2} \\
      2000 & ~~~~7 & {}   & {}              & ~\textbf{2.2} & ~0.2          & ~~22         & 1.8          & ~0.2          \\
      5000 & ~~~~7 & {}   & {}              & ~2.8          & ~\textbf{0.2} & ~~~7         & 2.8          & ~0.2          \\
      \bottomrule
    \end{tabular}
  \end{center}
\end{table}

\begin{table}[tp]
  \small
  \caption{Results for \texttt{hickerson-13}}
  \label{hickerson-13}
  \begin{center}
    \def~{\hphantom{0}}
    \begin{tabular}{cc*4{>{\sffamily}c}c*2{>{\sffamily}c}}
      \toprule
      & \multicolumn{5}{c}{Without irrationalization} 
      & \multicolumn{3}{c}{With irrationalization} \\
      \cmidrule{2-6} \cmidrule(l){7-9}
      & & \multicolumn{4}{c}{Time (s)} 
      & & \multicolumn{2}{c}{Time (s)} \\
      \cmidrule{3-6} \cmidrule{8-9}
      & & & & \multicolumn{2}{c}{LattE macchiato}
      & & \multicolumn{2}{c}{LattE macchiato} \\
      \cmidrule{5-6} \cmidrule{8-9}
      \smash[t]{\begin{tabular}[b]{c}Max.\\index\end{tabular}} 
      & Cones & \multicolumn{1}{c}{\smash[t]{\begin{tabular}[b]{@{}c@{}}LattE\\
            v\,1.2\end{tabular}} }
       & \multicolumn{1}{c}{\smash[t]{\begin{tabular}[b]{@{}c@{}}\texttt{barv.}\\
             v\,0.21\end{tabular}} }
      & \multicolumn{1}{c}{Poly} & \multicolumn{1}{c}{Exp} & Cones &
      \multicolumn{1}{c}{Poly} & \multicolumn{1}{c}{Exp} \\ 
      \midrule
      {~~~~1}  & 466\,540 & 793 & 589           & ~421          & 707          & 483\,507     & 479          & 770          \\  
      {~~~10}  & 272\,922 & {}  & {}            & ~\textbf{345} & 428          & ~55\,643     & 117          & 109          \\             
      {~~100}  & 142\,905 & {}  & {}            & ~489          & 249          & ~~9\,158     & \textbf{~83} & ~22          \\              
      {~~200}  & 122\,647 & {}  & {}            & ~625          & 222          & ~~6\,150     & ~93          & ~17          \\
      {~~300}  & ~98\,654 & {}  & {}            & ~903          & 199          & ~~4\,674     & 105          & ~14          \\
      {~~500}  & ~90\,888 & {}  & {}            & 1056          & 193          & ~~3\,381     & 137          & ~13 \\
      {~1000}  & ~73\,970 & {}  & {}            & 1648          & \textbf{190} & ~~2\,490     & 174          & ~\textbf{13} \\              
      {~2000}  & ~66\,954 & {}  & {}            & 2166          & 201          & ~~1\,857     & 237          & ~14          \\
      {~5000}  & ~49\,168 & {}  & {}            & 5040          & 286          & ~~1\,488     & 354          & ~18          \\
      {10000}  & ~43\,511 & {}  & {}            & 7278          & 370          & ~~1\,011     & 772          & ~34          \\              
      \bottomrule
    \end{tabular}
  \end{center}
\end{table}

\begin{table}[tp]
  \small
  \caption{Results for \texttt{hickerson-14}}
  \label{hickerson-14}
  \begin{center}
    \def~{\hphantom{0}}
    \begin{tabular}{cc*4{>{\sffamily}c}c*2{>{\sffamily}c}}
      \toprule
      & \multicolumn{5}{c}{Without irrationalization} 
      & \multicolumn{3}{c}{With irrationalization} \\
      \cmidrule{2-6} \cmidrule(l){7-9}
      & & \multicolumn{4}{c}{Time (s)} 
      & & \multicolumn{2}{c}{Time (s)} \\
      \cmidrule{3-6} \cmidrule{8-9}
      & & & & \multicolumn{2}{c}{LattE macchiato}
      & & \multicolumn{2}{c}{LattE macchiato} \\
      \cmidrule{5-6} \cmidrule{8-9}
      \smash[t]{\begin{tabular}[b]{c}Max.\\index\end{tabular}} 
      & Cones & \multicolumn{1}{c}{\smash[t]{\begin{tabular}[b]{@{}c@{}}LattE\\
            v\,1.2\end{tabular}} }
       & \multicolumn{1}{c}{\smash[t]{\begin{tabular}[b]{@{}c@{}}\texttt{barv.}\\
             v\,0.21\end{tabular}} }
      & \multicolumn{1}{c}{Poly} & \multicolumn{1}{c}{Exp} & Cones &
      \multicolumn{1}{c}{Poly} & \multicolumn{1}{c}{Exp} \\ 
      \midrule
      {~~~~1}  & 1\,682\,743 & 4\,017 & 15\,284 & ~2\,053        & ~3\,466          & 552\,065   & ~~\,792          & ~1\,244          \\
      {~~~10}  & 1\,027\,619 & {}     & {}      & ~\textbf{1736} & ~2\,177          & ~49\,632   & ~~\,168          & ~~\,143          \\       
      {~~100}  & ~\,455\,474 & {}     & {}      & ~2\,294        & ~1\,089          & ~~8\,470   & ~~\,\textbf{128} & ~~\,~29          \\        
      {~~200}  & ~\,406\,491 & {}     & {}      & ~2\,791        & ~~\,990          & ~~5\,554   & ~~\,157          & ~~\,~22          \\
      {~~300}  & ~\,328\,340 & {}     & {}      & ~4\,131        & ~~\,875          & ~~4\,332   & ~~\,187          & ~~\,~19          \\
      {~~500}  & ~\,303\,566 & {}     & {}      & ~4\,911        & ~~\,842          & ~~3\,464   & ~~\,235          & ~~\,~\textbf{18} \\
      {~1000}  & ~\,236\,626 & {}     & {}      & ~8\,229        & ~~\,\textbf{807} & ~~2\,384   & ~~\,337          & ~~\,~18          \\
      {~2000}  & ~\,195\,368 & {}     & {}      & 12\,122        & ~~\,817          & ~~1\,792   & ~~\,481          & ~~\,~21          \\
      {~5000}  & ~\,157\,496 & {}     & {}      & 22\,972        & ~1\,034          & ~~1\,276   & ~~\,723          & ~~\,~27          \\
      {10000}  & ~\,128\,372 & {}     & {}      & 31\,585        & ~1\,270          & ~~~\,956   & ~1\,095          & ~~\,~38          \\
      \bottomrule
    \end{tabular}
  \end{center}
\end{table}

\smallbreak
We can make the following observations:
\begin{enumerate}[\ i.]
\item By stopping Barvinok decomposition before the cones are unimodular, it
  is possible to significantly reduce the number of simplicial cones.  This
  effect is much stronger with irrational decomposition in the primal space than
  with decomposition in the dual space. 

\item The newly implemented exponential substitution has a computational
  overhead compared to the polynomial substitution that was implemented in
  LattE~1.2. 

\item However, when we compute with simplicial, non-unimodular cones, 
  the exponential substitution becomes much more efficient than the
  polynomial substitution.   Hence the break-even point between enumeration and
  decomposition is reached at a larger cone index.

  \noindent The reason is that the inner loops are shorter for the exponential
  substitution; essentially, only a sum of powers of scalar products needs to
  be evaluated in the formula~\eqref{eq:count-with-todd}.  This can be done
  very efficiently.  

\item The best results are obtained with the irrational primal decomposition
  down to an index of about 500 to 1000 and exponential substitution.

\end{enumerate} 

\subsection{Results for challenge problems}
In \autoref{tab:hickerson-table} we show the results for some larger test cases
related to \citet{hickerson-1991}.  We compare LattE 1.2 with our
implementation of irrational primal decomposition
(\autoref{algo:primal-irr-barvi}) with maximum index 500.  The computation
times are given in CPU seconds.  The computations with LattE 1.2 were done on
a PC Pentium M, 1.4\,GHz; the computations with LattE macchiato were
done on a slightly slower machine, a Sun Fire V890 with UltraSPARC-IV
processors, 1.2\,GHz.

\begin{table}[tb]
\centering
\small
\caption{Results for larger Hickerson problems}
\label{tab:hickerson-table}
\def~{\hphantom{0}}
\begin{tabular}{cccr>{\sffamily}rr>{\sffamily}r}
  \toprule
  & & & \multicolumn{2}{c}{LattE v\,1.2} & \multicolumn{2}{c}{LattE macchiato} \\
  \cmidrule{4-5}\cmidrule(l){6-7}
  $n$ & $d$ & \smash[t]{\begin{tabular}[b]{@{}c@{}}Lattice\\ points\end{tabular}}
  & Cones\hspace{1em}\mbox{} & \multicolumn{1}{c}{Time} 
  & Cones\hspace{1em}\mbox{} & \multicolumn{1}{c}{Time} \\
  \midrule 
  15   & ~7 & 20 & ~\,293\,000        & 10\,min 55\,s                      & 2\,000           & 22\,s                         \\
  16   & ~8 & 54 & 3\,922\,000        & 3\,h 35\,min \hphantom{00\,s}      & 19\,000          & 3\,min 56\,s                  \\
  17   & ~8 & 18 & {}                 & {}                                 & 2\,655\,000      & 7\,h 59\,min \hphantom{00\,s} \\
  18   & ~9 & 44 & 61\,500\,000       & 77\,h 00\,min \hphantom{00\,s}     & 200\,000         & 49\,min 12\,s                 \\
  20   & 10 & 74 & {}                 & {}                                 & 2\,742\,000      & 13\,h 05\,min \hphantom{00\,s} \\
  \bottomrule
\end{tabular}
\end{table}%

\begin{table}[tb]
\centering
\small
\caption{Results for cross polytopes}
\label{tab:cross-polytope}
\def~{\hphantom{0}}
\begin{tabular}{cc>{\sffamily}cc>{\sffamily}c}
  \toprule
  & \multicolumn{2}{c}{Without irrationalization} &
  \multicolumn{2}{c}{All-primal irrational} \\
  \cmidrule{2-3}\cmidrule(l){4-5}
  $d$ 
  & Cones\hspace{.5em}\mbox{} & \multicolumn{1}{c}{Time (s)} 
  & Cones\hspace{1em}\mbox{} & \multicolumn{1}{c}{Time (s)} \\
  \midrule 
  4 & ~~~\,384 & ~~~1.1  & & ~~\,~~0.9 \\
  5 & ~~3\,840 & ~~~6.5 & & ~~\,~~1.4\\
  6 & ~46\,264 & ~~91.7  & & ~~\,~~2.7 \\
  7 & 653\,824 & 1688.7 & & ~~\,~~5.5 \\
  8 &  &  & ~~1\,000 & ~~\,~12.3 \\
  9 & & & ~~2\,000 & ~~\,~29.6 \\
  10 & & & ~~5\,000 & ~~\,~74.8 \\
  11 & & & ~11\,000 & ~~\,189.1 \\
  12 & & & ~24\,000 & ~~\,483.0 \\
  13 & & & ~53\,000 & ~1\,231.2 \\
  14 & & & 114\,000 & ~3\,145.6 \\
  15 & & & 245\,000 & ~8\,180.9 \\
  \bottomrule
\end{tabular}
\end{table}%

Both the traditional Barvinok algorithm (\autoref{algo:dual-barvi} with $\ell=1$)
and the homogenized variant of Barvinok's algorithm \citep{latte2} do not work
well for cross polytopes. 
The reason is that triangulation is done in the dual space, so hypercubes need to
be triangulated.  We show the performance of the
traditional Barvinok algorithm in \autoref{tab:cross-polytope}.
We also show computational results for the all-primal
irrational algorithm (\autoref{algo:all-primal-irr-barvi} with $\ell=500$), using
exponential substitution.  The computation times are given in CPU seconds on a
Sun Fire V440 with UltraSPARC-IIIi processors, 1.6\,GHz. 

A challenge problem related to the paper by~\citet{beck-hosten:cyclotomic},
case $m=42$, could
be solved using the all-primal irrational decomposition algorithm
(\autoref{algo:all-primal-irr-barvi}) with
exponential substitution.  The method decomposed the polyhedron to a total of
1.1 million simplicial cones of index at most 500.  The computation took
66\,000 CPU seconds on a Sun Fire V440 with UltraSPARC-IIIi processors,
1.6\,GHz.  The problem could not be solved previously because the traditional 
algorithms first tried to triangulate the polar cones, which does not finish
within 17~days of computation.

\section*{Conclusions and future work}

The above computational results with our preliminary implementation have shown
that the proposed irrationalization techniques can speed up the Barvinok
algorithm by large factors. 

A further speed-up can be expected from a refined implementation.  For
example, the choice of the irrational shifting vector is based on worst-case
estimates.  It may be worthwhile to implement a randomized choice of the
shifting vector (within the stability cube), using shorter rational numbers
than those constructed in the paper.  The randomized choice, of course, would
not give the same guarantees as our deterministic construction.  However, it
is easy and efficient to check, during the decomposition, if the generated
cones are all irrational; when they are not, one could choose a new random
shifting vector (or resort to the one constructed in this paper) and restart
the computation.

\bibliographystyle{plainnat}
\bibliography{barvinok,iba-bib,weismantel}
\end{document}